\documentclass{amsart}
\usepackage{amssymb}
\usepackage{mathrsfs}
\usepackage{stmaryrd} 
\usepackage{bbm}
\usepackage{oldgerm}
\usepackage[english]{babel}
\usepackage[T1]{fontenc}
\usepackage[latin1]{inputenc}
\usepackage[all]{xy}
\usepackage{hyperref}
\usepackage{upref}
\usepackage{xcolor}
\usepackage{microtype}
\usepackage[Symbol]{upgreek}

\hypersetup{
    colorlinks,
    linkcolor={red!50!black},
    citecolor={blue!50!black},
    urlcolor={blue!80!black}
}

\newtheorem{teo}[subsection]{Theorem}
\newtheorem{prop}[subsection]{Proposition}
\newtheorem{cor}[subsection]{Corollary}

\theoremstyle{definition}

\newtheorem{defi}[subsection]{Definition}
\newtheorem{rema}[subsection]{Remark}

\numberwithin{equation}{subsection}

\mathchardef\mhyphen="2D

\DeclareMathSymbol{\mlq}{\mathord}{operators}{``}
\DeclareMathSymbol{\mrq}{\mathord}{operators}{`'}

\newcommand{\mH}{{\mathbb H}}

\newcommand{\mQ}{{\mathbb Q}}

\newcommand{\mN}{{\mathbb N}}

\newcommand{\mX}{{\mathbb X}}

\newcommand{\mZ}{{\mathbb Z}}

\newcommand{\mG}{{\mathbb G}}
\newcommand{\mK}{{\mathbb K}}

\newcommand{\mV}{{\mathbb V}}

\newcommand{\bB}{{\bf B}}

\newcommand{\bD}{{\bf D}}

\newcommand{\bIndMod}{{\bf Ind\mhyphen Mod}}

\newcommand{\Et}{{\bf \acute{E}t}}

\newcommand{\bRep}{{\bf Rep}}

\newcommand{\bMod}{{\bf Mod}}

\newcommand{\colim}{{\underset{\longrightarrow}{\lim}}}
\newcommand{\indcolim}{{\mlq\mlq\colim \mrq\mrq}}

\newcommand{\et}{{\rm \acute{e}t}}
\newcommand{\fet}{{\rm f\acute{e}t}}

\newcommand{\sm}{{\rm sm}}

\newcommand{\zar}{{\rm zar}}

\newcommand{\coh}{{\rm coh}}

\newcommand{\HT}{{\rm HT}}
\newcommand{\Dolb}{{\rm Dolb}}
\newcommand{\sol}{{\rm sol}}
\newcommand{\sDolb}{{\rm sDolb}}

\newcommand{\qsolnilp}{{\rm qsolnilp}}

\newcommand{\rf}{{\rm f}}

\newcommand{\Spec}{{\rm Spec}}

\newcommand{\Spf}{{\rm Spf}}

\newcommand{\ob}{{\rm Ob}}

\newcommand{\tot}{{\rm tot}}

\newcommand{\cont}{{\rm cont}}

\newcommand{\id}{{\rm id}}

\newcommand{\Sym}{{\rm Sym}}

\newcommand{\Hom}{{\rm Hom}}

\newcommand{\bHM}{{\bf HM}}

\newcommand{\rC}{{\rm C}}

\newcommand{\rE}{{\rm E}}

\newcommand{\rH}{{\rm H}}
\newcommand{\rI}{{\rm I}}

\newcommand{\rR}{{\rm R}}

\newcommand{\rW}{{\rm W}}

\newcommand{\oK}{{\overline{K}}}

\newcommand{\oR}{{\overline{R}}}
\newcommand{\oS}{{\overline{S}}}

\newcommand{\oU}{{\overline{U}}}

\newcommand{\oX}{{\overline{X}}}

\newcommand{\ox}{{\overline{x}}}
\newcommand{\oy}{{\overline{y}}}

\newcommand{\oeta}{{\overline{\eta}}}

\newcommand{\ocB}{{\overline{\cB}}}

\newcommand{\uX}{{\underline{X}}}

\newcommand{\ucH}{{\underline{\cH}}}

\newcommand{\upi}{{\underline{\pi}}}

\newcommand{\umK}{{\underline{\mK}}}

\newcommand{\hA}{{\widehat{A}}}

\newcommand{\hRun}{{\widehat{R_1}}}

\newcommand{\halpha}{\widehat{\alpha}}

\newcommand{\hupsigma}{\widehat{\upsigma}}
\newcommand{\hupbeta}{\widehat{\upbeta}}

\newcommand{\vupsigma}{\vec{\upsigma}}

\newcommand{\cA}{{\mathscr A}}
\newcommand{\cB}{{\mathscr B}}
\newcommand{\cC}{{\mathscr C}}

\newcommand{\cF}{{\mathscr F}}
\newcommand{\cG}{{\mathscr G}}

\newcommand{\cL}{{\mathscr L}}

\newcommand{\co}{{\mathscr O}}

\newcommand{\cS}{{\mathscr S}}

\newcommand{\cH}{{\mathscr H}}
\newcommand{\cM}{{\mathscr M}}
\newcommand{\cN}{{\mathscr N}}

\newcommand{\cV}{{\mathscr V}}

\newcommand{\fX}{{\mathfrak X}}

\newcommand{\fd}{{\mathfrak d}}
\newcommand{\fgg}{{\mathfrak g}}
\newcommand{\fm}{{\mathfrak m}}

\newcommand{\hoR}{{\widehat{\oR}}}

\newcommand{\hcC}{{\widehat{\cC}}}

\newcommand{\bvg}{{\breve{g}}}

\newcommand{\bvocB}{{\breve{\ocB}}}

\newcommand{\bvcC}{{{\breve{\cC}}}}

\newcommand{\bvcF}{{\breve{\cF}}}

\newcommand{\bvpsi}{{\breve{\psi}}}

\newcommand{\bvuptheta}{{\breve{\uptheta}}}

\newcommand{\bvmZ}{{\breve{\mZ}}}

\newcommand{\tE}{{\widetilde{E}}}

\newcommand{\tG}{{\widetilde{G}}}

\newcommand{\tS}{{\widetilde{S}}}

\newcommand{\tU}{{\widetilde{U}}}

\newcommand{\tX}{{\widetilde{X}}}

\newcommand{\tg}{{\widetilde{g}}}

\newcommand{\tx}{{\widetilde{x}}}

\newcommand{\txi}{{\widetilde{\xi}}}

\newcommand{\tmX}{{\widetilde{\mX}}}

\begin{document}

\title[The p-adic Simpson Correspondence II]{The $p$-adic Simpson Correspondence II: Functoriality by proper direct image and Hodge-Tate local systems -- \\ an overview}
\author{Ahmed Abbes and Michel Gros}
\address{A.A. Laboratoire Alexander Grothendieck, UMR 9009 du CNRS, 
Institut des Hautes \'Etudes Scientifiques, 35 route de Chartres, 91440 Bures-sur-Yvette, France}
\address{M.G. CNRS UMR 6625, IRMAR, Université de Rennes 1,
Campus de Beaulieu, 35042 Rennes cedex, France}
\email{abbes@ihes.fr}
\email{michel.gros@univ-rennes1.fr}

\begin{abstract}
Faltings initiated in 2005 a $p$-adic analogue of the (complex) Simpson correspondence whose construction has been taken up 
by various authors, according to several approaches.
Following the one we initiated in \cite{agt}, we present an overview of a new monograph \cite{ag2} developing new features of 
the $p$-adic Simpson correspondence, inspired by our construction of the relative Hodge-Tate spectral sequence \cite{ag1}. 
First, we address the connection to Hodge-Tate local systems. Second, we establish the functoriality of the $p$-adic Simpson correspondence by proper direct image. 
Along the way, we expand the scope of our original construction. 
\end{abstract}

\maketitle

\setcounter{tocdepth}{1}
\tableofcontents

\section{Introduction}\label{intro}

\subsection{}\label{intro1}
Initiated by Faltings \cite{faltings3} and developed following different approaches including those by Tsuji and by the authors \cite{agt}, 
the $p$-adic Simpson correspondence provides an equivalence of categories between certain {\em $p$-adic étale local systems} 
over a smooth algebraic variety over a local field and certain {\em Higgs bundles}. 
The key idea behind its construction comes from Faltings' approach in $p$-adic Hodge theory, more precisely, from his strategy for the computation 
of the cohomology of a $p$-adic étale local system. 

\subsection{}\label{intro3}
Let $K$ be a complete discrete valuation field of characteristic $0$, with {\em algebraically closed} residue field $k$ of characteristic $p>0$, 
$\co_K$ the valuation ring of $K$, $\oK$ an algebraic closure of $K$, $\co_\oK$ the integral closure of $\co_K$ in $\oK$,
$G_K$ the Galois group of $\oK$ over $K$, $\co_C$ the $p$-adic completion of $\co_\oK$, $\fm_C$ its maximal ideal, $C$ its field of fractions.
We set $S=\Spec(\co_K)$ and $\oS=\Spec(\co_\oK)$ and we denote by $s$ (resp.  $\eta$, resp. $\oeta$) 
the closed point of $S$ (resp. generic point of $S$, resp. generic point of $\oS$).

Let $X$ be a proper smooth $S$-scheme and let $L$ be a locally constant constructible sheaf of 
$\mZ/p^n\mZ$-modules of $X_{\oeta,\et}$ for an integer $n\geq 0$. 
To compute the cohomology of $L$, Faltings introduced a ringed topos $(\tE,\ocB)$ equipped with two morphisms of topos 
\begin{equation}
\xymatrix{
X_{\oeta,\et}\ar[r]^-(0.5){\psi}&{\tE}\ar[r]^-(0.5){\sigma}&X_\et}.
\end{equation}
For any integer $j\geq 1$, we have $\rR^j\psi_*(L)=0$. In particular, for any $i\geq 0$, we have a canonical isomorphism
\begin{equation}\label{intro3a}
\rH^i(X_{\oeta,\et},L)\stackrel{\sim}{\rightarrow}\rH^i(\tE,\psi_*(L)).
\end{equation} 
Using Artin-Schreier theory, Faltings proved a refinement of this result, namely that the canonical morphism 
\begin{equation}\label{intro3b} 
\rH^i(X_{\oeta,\et},L)\otimes_{\mZ_p}\co_C\rightarrow \rH^i(\tE,\psi_*(L)\otimes_{\mZ_p}\ocB)
\end{equation}
is an {\em almost isomorphism}, {\em i.e.}, its kernel and cokernel are annihilated by $\fm_C$.  
Setting $\cM=\psi_*(L)\otimes_{\mZ_p}\ocB$, we have a canonical isomorphism 
\begin{equation}\label{intro3c}
\rR\Gamma(\tE,\cM) \stackrel{\sim}{\rightarrow}\rR\Gamma(X_\et,\rR\sigma_*(\cM)). 
\end{equation}
This computation extends naturally to $\mQ_p$-local systems by passing to the limit on $n$ and inverting $p$. 
For certain $\mQ_p$-local systems $L$, the corresponding complex $\rR\sigma_*(\cM)$ 
is the Dolbeault complex of a Higgs bundle canonically associated to $L$ by the {\em $p$-adic Simpson correspondence} 
(see \cite{agt} II.2.8 and \cite{ag2} 2.8.1 for the terminology on Higgs modules). The Cartan-Leray spectral sequence for $\sigma$ \eqref{intro3c} 
and the almost isomorphism \eqref{intro3b} lead then to a generalization of the Hodge-Tate spectral sequence \cite{ag1}. 

\subsection{}
To describe our construction of the $p$-adic Simpson correspondence, 
we consider the ring $\bvocB=(\ocB/p^n\ocB)_{n\geq 0}$ of the topos $\tE^{\mN^\circ}$ of projective systems of objects of $\tE$, 
indexed by the ordered set $\mN$. We work in the {\em category of $\bvocB$-modules up to isogeny}, that we call the category of $\bvocB_\mQ$-modules. 
It is a counterpart on the side of Faltings topos of the category of $\mQ_p$-local systems of $X_{\oeta,\et}$.

\subsection{}\label{intro4}
We built in \cite{agt} a correspondence between certain $\bvocB_\mQ$-modules and certain Higgs bundles, through the classical scheme of Fontaine correspondences
involving a period ring of $\tE$ that we call the {\em Higgs-Tate algebra}. It is an integral model of Hyodo's ring $B_\HT$, that we construct using deformation theory, 
inspired by Faltings' original approach and the work of Ogus-Vologodsky on the Cartier transform in characteristic $p$ \cite{ov}. 
A ``weak'' $p$-adic completion of the Higgs-Tate algebra has nice cohomological properties that lead to an equivalence between the categories of admissible objects, 
namely {\em Dolbeault $\bvocB_\mQ$-modules} and {\em solvable Higgs bundles}. 

\subsection{}\label{intro5}
We review in the new monograph \cite{ag2} the construction of this correspondence enlarging along the way its scope, to prove the following important new features: 
\begin{itemize}
\item[(i)] We characterize the {\em Hodge-Tate} $\bvocB_\mQ$-modules among the Dolbeault $\bvocB_\mQ$-modules, namely, those whose associated Higgs bundle is nilpotent. 
\item[(ii)] We prove the functoriality of the $p$-adic Simpson correspondence by proper direct image, which leads to a generalization of the relative Hodge-Tate spectral sequence \cite{ag1}.
\end{itemize}  

Another new feature worth mentioning is the cohomological descent of Dolbeault modules in the small affine case (\cite{ag2} 4.8.16) 
which leads to a description of the category of Dolbeault modules in terms of their global sections, and allows to compare the global and the local theories (\cite{ag2} 4.8.31 and 4.8.32).

\subsection{}\label{intro6}
The $p$-adic Simpson correspondence requires the existence of a smooth deformation of the scheme $X\otimes_{\co_K}\co_C$ 
over Fontaine's universal $p$-adic $\rW(k)$-infinitesimal thickening $\cA_2(\co_\oK)$ (see \cite{agt} (II.9.3.5), \cite{ag2} 2.3.5 and 2.3.8), and it depends on the choice of such a deformation, 
while the theory of Hodge-Tate $\bvocB_\mQ$-modules does not depend on it. This apparent paradox is resolved by introducing a new version of our construction 
of the $p$-adic Simpson correspondence that depends on a smooth deformation of the scheme $X\otimes_{\co_K}\co_C$ 
over a logarithmic version $\cA^*_2(\co_\oK/\co_K)$ of Fontaine's universal $p$-adic $\co_K$-infinitesimal thickening $\cA_2(\co_\oK/\co_K)$ (cf. \cite{ag2} (2.3.2.4) and (2.4.5.1)). 
We refer to this context as the {\em relative case} and to the former context as the {\em absolute case}.
We deal with both cases simultaneously, each having its advantages and disadvantages.
Since $\cA^*_2(\co_\oK/\co_K)$ is naturally an $\co_K$-algebra, $X\otimes_{\co_K}\co_C$ has a canonical smooth $\cA^*_2(\co_\oK/\co_K)$-deformation, namely the base change of $X$. 
The theory in the relative case has its limits, however. Indeed, we know in the affine case that Dolbeault modules are small (\cite{ag2} 3.4.30); 
roughly speaking, they are ``trivial'' modulo a prescribed power of $p$. This power of $p$ is independent of $K$ only in the absolute case. Its dependence on $K$ 
limits drastically the scope of the theory in the relative case, in particular, if we would like to extend the correspondence to all $\bvocB_\mQ$-modules by descent.  
However, Hodge-Tate $\bvocB_\mQ$-modules are Dolbeault both in the absolute and the relative cases, 
which makes our definition of Hodge-Tate $\bvocB_\mQ$-modules independent of any deformation. 

\subsection{}\label{intro7}
Compared to \cite{agt}, we also extend the theory in \cite{ag2} from the category of $\bvocB_\mQ$-modules to the larger category of {\em ind-$\bvocB$-modules}. 
The latter admits filtered inductive limits and has better properties (\cite{ks2} and \cite{ag2} 2.6 and 2.7).  
This enlargement which seems at first sight technical, turns out to be quite useful and is necessary for the study of the functoriality of the $p$-adic 
Simpson correspondence by proper direct image.

\subsection{}\label{intro8}
We give below a detailed overview of the content of our  new monograph \cite{ag2}. 
We treat in {\em loc. cit.} schemes with toric singularities using logarithmic geometry, 
but for simplicity, we stick in this overview to the smooth case.  

\subsection*{Acknowledgments} We would like to warmly thank T. Tsuji for the many advices he very generously gave us throughout this work.

\section{Faltings topos}\label{ft}

\subsection{}\label{ft1}
Let $X$ be a smooth $S$-scheme,  $E$ the category of morphisms $(V\rightarrow U)$
above the canonical morphism $X_\oeta\rightarrow X$, that is, commutative diagrams 
\begin{equation}\label{ft1a}
\xymatrix{V\ar[r]\ar[d]&U\ar[d]\\
X_\oeta\ar[r]&X}
\end{equation}
such that $U$ is étale over $X$ and the canonical morphism $V\rightarrow U_\oeta$ is {\em finite étale}. 
It is useful to consider the category $E$ as fibred by the functor
\begin{equation}\label{ft1b}
\pi\colon E\rightarrow \Et_{/X}, \ \ \ (V\rightarrow U)\mapsto U,
\end{equation}
over the étale site of $X$. 

The fiber of $\pi$ above an object $U$ of $\Et_{/X}$ is canonically equivalent to the category $\Et_{\rf/U_\oeta}$ of finite étale morphisms 
over $U_\oeta$. We equip it with the étale topology and denote by $U_{\oeta,\fet}$ the associated topos. 
If $U_\oeta$ is connected and if $\oy$ is a geometric point of $U_\oeta$, then the topos $U_{\oeta,\fet}$ is equivalent to
the classifying topos of the profinite group $\pi_1(U_\oeta,\oy)$, {\em i.e.}, the category of discrete sets equipped with a continuous left action
of $\pi_1(U_\oeta,\oy)$.

\subsection{}\label{ft2}
We equip $E$ with the {\em covanishing} topology (\cite{agt} VI.10.1), that is the topology generated by coverings 
$\{(V_i\rightarrow U_i)\rightarrow (V\rightarrow U)\}_{i\in I}$
of the following two types~:
\begin{itemize}
\item[(v)] $U_i=U$ for all $i\in I$ and $(V_i\rightarrow V)_{i\in I}$ is a covering;
\item[(c)] $(U_i\rightarrow U)_{i\in I}$ is a covering and $V_i=V\times_UU_i$ for all $i\in I$. 
\end{itemize}

We denote by $\tE$ the topos of sheaves of sets on $E$. 

To give a sheaf $F$ on $E$ is equivalent to give:
\begin{itemize}
\item[(i)] for any object $U$ of $\Et_{/X}$, a sheaf $F_U$ of $U_{\oeta,\fet}$, namely the restriction of $F$ to the fiber
of $\pi$ above $U$;
\item[(ii)] for any morphism $f\colon U'\rightarrow U$ of $\Et_{/X}$, a morphism $\gamma_f\colon F_U\rightarrow f_{\oeta*}(F_{U'})$. 
\end{itemize}
These data should satisfy a cocycle condition (for the composition of morphisms) and  
a gluing condition (for coverings of $\Et_{/X}$). We write $F=\{U \mapsto F_U\}$ and think of it as a sheaf on $\Et_{/X}$ with values in finite étale topos.

\subsection{}\label{ft20}
Any specialization map $\oy\rightsquigarrow \ox$ from a geometric point $\oy$ of $X_\oeta$ to a geometric point $\ox$ of $X$, 
determines a point of $\tE$ that we denote by $\rho(\oy\rightsquigarrow \ox)$. {\em The collection of these points of $\tE$ is conservative.} 

\subsection{}\label{ft3}
There are three morphisms of topos (\cite{agt} VI.10.6 and VI.10.7)
\begin{equation}\label{ft3a}
\xymatrix{
X_{\oeta,\et}\ar[r]^-(0.5){\psi}&{\tE}\ar[r]^-(0.5){\sigma}\ar[d]^{\beta}&X_\et\\
&X_{\oeta,\fet}&}
\end{equation}
defined by 
\begin{eqnarray}
(V\mapsto U)\in \ob(E) &\mapsto& \psi^*(V\rightarrow U)=V,\label{ft3d}\\
U\in \ob(\Et_{/X}) &\mapsto& \sigma^*(U)=(U_\oeta\rightarrow U)^a,\label{ft3b}\\
V\in \ob(\Et_{\rf/X_\oeta}) &\mapsto& \beta^*(V)=(V\rightarrow X)^a.\label{ft3c}
\end{eqnarray}

\subsection{}\label{ft30}
The higher direct images of $\sigma$ sheafify Galois cohomology (\cite{agt} VI.10.40): if $F=\{U\mapsto F_U\}$ is an abelian group of $\tE$,  
for each integer $i\geq 0$, $\rR^i\sigma_*(F)$ is canonically isomorphic to the sheaf associated to the presheaf 
\begin{equation}\label{ft30a}
U\mapsto \rH^i(U_{\oeta,\fet}, F_U).
\end{equation}

\begin{prop}\label{ft4}
For any locally constant constructible torsion abelian sheaf $F$ of $X_{\oeta,\et}$, we have $\rR^i\psi_*(F)=0$ for any $i\geq 1$.
\end{prop}

This statement is a consequence of the fact that for any geometric point $\ox$ of $X$ over $s$, denoting by $\uX$ the strict
localization of $X$ at $\ox$, $\uX_\oeta$ is a $K(\pi,1)$ scheme. This property was proved by Faltings (\cite{faltings1} Lemma 2.3 page 281), 
generalizing results of Artin (\cite{sga4} XI). It was generalized further by Achinger in the log-smooth case \cite{achinger}.

\subsection{}\label{ft5}
For any object $(V\rightarrow U)$ of $E$, we denote by $\oU^V$ the integral closure of $\oU=U\times_S\oS$ in $V$ and we set
\begin{equation}\label{ft5a}
\ocB(V\rightarrow U)=\Gamma(\oU^V,\co_{\oU^V}).
\end{equation} 
The presheaf on $E$ defined above is in fact a sheaf. Like any sheaf, we can write $\ocB=\{U\mapsto \ocB_U\}$.

Let $U=\Spec(R)$ be an étale $X$-scheme, $\oy$ a geometric point of $U_\oeta$. The stalk $\ocB_{U,\oy}$ can be described as follows.   
We denote by $(V_i)_{i\in I}$ the universal cover of $U_\oeta$ at $\oy$.
For each $i\in I$, let $U_i=\Spec(R_i)$ be the normalization of $\oU$ in $V_i$
\begin{equation}\label{ht5b}
\xymatrix{
V_i\ar[r]\ar[d]&U_i\ar[d]\\
U_\oeta\ar[r]&\oU}
\end{equation}
Then, the stalk $\ocB_{U,\oy}$ is isomorphic to the following $\co_\oK$-representation of $\pi_1(U_\oeta,\oy)$:
\begin{equation}\label{ht5c}
\oR=\underset{\underset{i\in I}{\longrightarrow}}{\lim} \ R_i. 
\end{equation}

\vspace{2mm}

Using Artin-Schreier theory, Faltings proved the following refinement of \ref{ft4}:

\begin{teo}[\cite{faltings2}, \cite{ag1} 4.8.13]\label{ft6}
For any locally constant constructible sheaf of $(\mZ/p^n\mZ)$-modules $F$ of $X_{\oeta,\et}$, the canonical morphism 
\begin{equation}\label{grfmc2a} 
\rH^i(X_{\oeta,\et},F)\otimes_{\mZ_p}\co_C\rightarrow \rH^i(\tE,\psi_*(F)\otimes_{\mZ_p}\ocB)
\end{equation}
is an {\em almost isomorphism}, i.e., its kernel and cokernel are annihilated by $\fm_C$. 
\end{teo}

Faltings derived all comparaison theorems between $p$-adic étale cohomology 
and other $p$-adic cohomologies from this main $p$-adic comparaison theorem.

\section{Local theory. The torsor of deformations} \label{tordef}

\subsection{}
First we review the local variant of the $p$-adic Simpson correspondence for {\em small} affine $S$-schemes.
Our approach uses a period ring, the {\em Higgs-Tate algebra}, that we introduced in \cite{agt}.

\subsection{}\label{tordef1}
Recall (\cite{agt} II.9, \cite{tsuji1} 1.1) that Fontaine associated functorially to any $\mZ_{(p)}$-algebra $A$ such that $A/pA\not=0$, the ring  
\begin{equation}\label{tordef1a}
A^\flat=\underset{\underset{\mN}{\longleftarrow}}{\lim}\ A/pA,
\end{equation} 
where the transition morphisms are the absolute Frobenius morphisms of $A/pA$, and the ring homomorphism
\begin{equation}\label{tordef1b}
\theta\colon \rW(A^\flat)\rightarrow \hA,
\end{equation}
from the Witt vectors of $A^\flat$ to the $p$-adic completion of $A$, 
defined for $x=(x_0,x_1,\dots)\in \rW(A^\flat)$ by 
\begin{equation}
\theta(x)=\underset{m\rightarrow +\infty}{\lim}\ (\tx_{0m}^{p^m}+p\tx_{1m}^{p^{m-1}}+\dots+p^m\tx_{mm}),
\end{equation}
where for each $n\geq 0$, we write $x_n=(x_{nm})_{m\geq 0}\in A^\flat$ and for $x\in A/pA$, $\tx$ denotes a lifting in $A$. 

The ring $A^\flat$ is perfect of characteristic $p$, and the homomorphism $\theta$ is surjective if the absolute Frobenius morphism of $A/pA$ is surjective.

\subsection{}\label{tordef2}
The ring $\co_{\oK^\flat}=(\co_\oK)^\flat$ is a complete non-discrete valuation ring of height $1$. We denote by $\oK^\flat$ its fractions field.  
We choose a sequence $(p_n)_{n\geq 0}$ of elements of $\co_\oK$ such that $p_0=p$ and $p_{n+1}^p=p_n$ for all $n\geq 0$. 
We denote by $\varpi$ the element of $\co_{\oK^\flat}$ defined by $(p_n)$ and we set 
\begin{equation}\label{tordef2a}
\xi=[\varpi]-p\in \rW(\co_{\oK^\flat}).
\end{equation}
It is a generator of the kernel of $\theta$ \eqref{tordef1b}. We set 
\begin{equation}\label{tordef2c}
\cA_2(\co_\oK)=\rW(\co_{\oK^\flat})/\ker(\theta)^2.
\end{equation}
Then, we have an exact sequence
\begin{equation}\label{tordef2d}
0\longrightarrow \co_C\stackrel{\cdot \xi}{\longrightarrow} \cA_2(\co_\oK)
\stackrel{\theta}{\longrightarrow} \co_C \longrightarrow 0.
\end{equation}

We have a canonical homomorphism $\mZ_p(1)\rightarrow \co_{\oK^\flat}^\times$. 
For any $\zeta\in \mZ_p(1)$, we have $\theta([\zeta]-1)=0$. We deduce a group homomorphism
\begin{equation}\label{tordef2e}
\mZ_p(1)\rightarrow \cA_2(\co_{\oK}),\ \ \ 
\zeta\mapsto\log([\zeta])=[\zeta]-1,
\end{equation}
whose image generates the ideal $p^{\frac{1}{p-1}}\xi\co_C$ de $\cA_2(\co_{\oK})$. It induces an $\co_C$-linear isomorphism (\cite{agt} II.9.18)
\begin{equation}\label{tordef2f}
\co_C(1)\stackrel{\sim}{\rightarrow} p^{\frac{1}{p-1}}\xi \co_C.
\end{equation}

\subsection{}\label{tordef20}
We introduce a  logarithmic relative version of the extension $\cA_2(\co_\oK)$ \eqref{tordef2d} over $\co_K$. 
We fix a uniformizer $\pi$ of $\co_K$ and a sequence $(\pi_n)_{n\geq 0}$ 
of elements of $\co_\oK$ such that $\pi_0=\pi$ and $\pi_{n+1}^p=\pi_n$ (for all $n\geq 0$) and let $\upi$ be the associated element of $\co_{\oK^\flat}$. 
We set
\begin{equation}\label{tordef20a}
\rW_{\co_K}(\co_{\oK^\flat})=\rW(\co_{\oK^\flat})\otimes_{\rW(k)}\co_K. 
\end{equation}
We denote by $\rW^{\ast}_{\co_K}(\co_{\oK^\flat})$ the sub-$\rW_{\co_K}(\co_{\oK^\flat})$-algebra of $\rW_K(\co_{\oK^\flat})=\rW(\co_{\oK^\flat})\otimes_{\rW(k)}K$ 
generated by $[\upi]/\pi$ and we set 
\begin{equation}\label{tordef20b}
\xi^{\ast}_\pi=\frac{[\upi]}{\pi}-1\in \rW^{\ast}_{\co_K}(\co_{\oK^\flat}).
\end{equation}
It is a generator of the kernel of the homomorphism $\theta^{\ast}_{\co_K}\colon \rW^{\ast}_{\co_K}(\co_{\oK^\flat})\rightarrow \co_C$ induced by $\theta$ \eqref{tordef1b}.
Observe that this algebra depends on $(\pi_n)_{n\geq 0}$. We set 
\begin{equation}\label{tordef20c}
\cA^{\ast}_2(\co_\oK/\co_K)=\rW^{\ast}_{\co_K}(\co_{\oK^\flat})/(\xi^{\ast}_\pi)^2 \rW^{\ast}_{\co_K}(\co_{\oK^\flat}).
\end{equation}
Then, we have an exact sequence 
\begin{equation}\label{tordef20d}
\xymatrix{
0\ar[r]&{\co_C}\ar[r]^-(0.5){\cdot \xi^{\ast}_\pi}&{\cA^{\ast}_2(\co_\oK/\co_K)}\ar[r]^-(0.5){\theta^{\ast}_{\co_K}}&{\co_C}\ar[r]& 0}.
\end{equation}

Denoting by $K_0$ the field of fractions of $\rW(k)$ and by $\fd$ the different of $K/K_0$, 
the canonical homomorphism $\cA_2(\co_\oK)\rightarrow \cA^{\ast}_2(\co_\oK/\co_K)$ induces an $\co_C$-linear isomorphism 
\begin{equation}\label{tordef20e}
\xi\co_C\stackrel{\sim}{\rightarrow}\pi\fd\xi^{\ast}_\pi \co_C.
\end{equation}

\subsection{}\label{tordef3}
Let $X=\Spec(R)$ be an affine smooth $S$-scheme which is {\em small}  in the sense of Faltings ({\em i.e.} it admits an étale $S$-morphism 
$X\rightarrow \mG_{m,S}^d=\Spec(\co_K[T_1^{\pm1},\dots, T_d^{\pm1}])$, for an integer $d\geq 0$),  
such that $X_s\not=\emptyset$. We fix  a geometric point $\oy$ of $X_\oeta$ and denote by $\oX^\star$ (resp. $X_\oeta^\star$) the connected component de
$\oX=X\times_S\oS$ (resp. $X_\oeta$) containing the image of $\oy$; we have  $X_\oeta^\star=\oX^\star\times_\oS\oeta$. 
We set $\Delta=\pi_1(X_\oeta^\star,\oy)$, $\oR=\ocB_{X,\oy}$ \eqref{ht5c} and 
\begin{equation}
R_1=\Gamma(\oX^\star,\co_{\oX}).
\end{equation}
Let $\hRun$ (resp. $\hoR$) be the $p$-adic completion of $R_1$ (resp. $\oR$).  Recall that $\Delta$ acts naturally on $\oR$ by ring homomorphisms. 
We denote by $\bRep_{\hoR}(\Delta)$ the category of $\hoR$-representations of $\Delta$ (cf. \cite{ag2} 2.1.2).

Applying the constructions of \ref{tordef1}--\ref{tordef20} to the algebra $\oR$ (\cite{ag2} 3.2.8 and 3.2.10), we set 
\begin{eqnarray}
\cA_2(\oR)=\rW(\oR^\flat)/\xi^2\rW(\oR^\flat),\label{tordef3b}\\
\cA^*_2(\oR/\co_K)=\rW^*_{\co_K}(\oR^\flat)/(\xi^{\ast}_\pi)^2 \rW_{\co_K}(\oR^\flat).\label{tordef3bb}
\end{eqnarray}
We have exact sequences
\begin{equation}\label{tordef3c}
0\longrightarrow \hoR\stackrel{\cdot \xi}{\longrightarrow} \cA_2(\oR)\longrightarrow \hoR \longrightarrow 0,
\end{equation}
\begin{equation}\label{tordef3cc}
0\longrightarrow \hoR\stackrel{\cdot \xi^{\ast}_\pi}{\longrightarrow} \cA^*_2(\oR/\co_K)\longrightarrow \hoR \longrightarrow 0.
\end{equation}

\subsection{}\label{tordef21}
Let $\tS$ be one of the schemes $\Spec(\cA_2(\co_\oK))$ or  $\Spec(\cA^{\ast}_2(\co_\oK/\co_K))$; the first case will be called {\em absolute} and the second {\em relative}. 
We denote by  
\begin{equation}\label{tordef21a}
i_S\colon \Spec(\co_C)\rightarrow \tS
\end{equation}
the closed immersion, defined by the ideal of square zero generated by $\txi=\xi$ in the absolute case and by $\txi=\xi^{\ast}_\pi$ in the relative case. 

Note that in the relative case, $\tS$ is naturally an $S$-scheme.

\subsection{}\label{tordef210}
Depending on whether we are in the absolute or relative case, we denote by $\tmX$ 
one of the schemes $\Spec(\cA_2(\oR))$ or $\Spec(\cA^*_2(\oR/\co_K))$. There is a canonical closed immersion 
\begin{equation}\label{tordef21b}
i_X\colon \Spec(\hoR)\rightarrow \tmX
\end{equation}
above the closed immersion $i_S$, defined by the ideal of $\co_\tmX$ generated in $\txi$. 

For any $\hRun$-algebra $A$, we consider Higgs $A$-modules with coefficients in $\txi^{-1}\Omega^1_{R/\co_K}\otimes_RA$ (cf. \cite{ag2} 2.5.1). 
We say abusively that they have coefficients in $\txi^{-1}\Omega^1_{R/\co_K}$. 
The category of these modules will be denoted by $\bHM(A,\txi^{-1}\Omega^1_{R/\co_K})$. 

\subsection{}\label{tordef4}
The $p$-adic Simpson correspondence depends on the choice of a smooth $\tS$-deformation $\tX$ of $X\otimes_{\co_K}\co_C$, 
\begin{equation}\label{tordef4a}
\xymatrix{
{X\otimes_{\co_K}\co_C}\ar[r]\ar[d]\ar@{}[rd]|\Box&{\tX}\ar[d]\\
{\Spec(\co_C)}\ar[r]&{\tS}}
\end{equation}
Since $X$ is affine, such a deformation always exists, and is unique up to a non-unique isomorphism. We fix one in the following. 

Let $U$ be an open subscheme of $\Spec(\hoR)$ and $\tU$ the open subscheme of
$\tmX$ defined by $U$ \eqref{tordef21b}.  We denote by $\cL(U)$ the set of morphisms represented 
by dotted arrows that complete the diagram 
\begin{equation}\label{tordef4b}
\xymatrix{
{U}\ar[d]\ar[r]&{\tU}\ar@{.>}[d]\ar@/^2pc/[dd]\\
{X\otimes_{\co_K}\co_C}\ar[r]\ar[d]\ar@{}[rd]|\Box&{\tX}\ar[d]\\
{\Spec(\co_C)}\ar[r]&{\tS}}
\end{equation}
in such a way that it remains commutative. The functor $U\mapsto \cL(U)$ is a torsor for the Zariski topology of $\Spec(\hoR)$ 
under the $\hoR$-module $\Hom_{R}(\Omega^1_{R/\co_K},\txi\hoR)$ (cf. \cite{ag2} 3.2.14). 
Such a torsor is very easy to describe. Let $\cF$ be the $\hoR$-module of affine functions on $\cL$ (\cite{agt} II.4.9). 
The latter fits into a canonical exact sequence
\begin{equation}\label{tordef4c}
0\rightarrow \hoR\rightarrow \cF\rightarrow \txi^{-1}\Omega^1_{R/\co_K} \otimes_R \hoR\rightarrow 0.
\end{equation} 
We consider the $\hoR$-algebra
\begin{equation}\label{tordef4d}
\cC=\underset{\underset{n\geq 0}{\longrightarrow}}\lim\ \Sym^n_{\hoR}(\cF),
\end{equation}
where the transition morphisms are defined by mapping $x_1\otimes\dots \otimes x_n$ to $1\otimes x_1\otimes\dots \otimes x_n$. 
Then, $\cL$ is represented by $\Spec(\cC)$ (\cite{agt} II.4.10). 

The natural action of $\Delta$ on $\oR$ induces an action on the scheme $\tmX$, and hence an $\hoR$-semi-linear action on $\cF$, 
such that the morphisms of \eqref{tordef4c} are $\Delta$-equivariant (cf. \cite{ag2} 3.2.15). 
We deduce an action of $\Delta$ on $\cC$ by ring homomorphisms. These actions are continuous for the $p$-adic topology (\cite{agt} II.12.4). 
We call the $\hoR$-algebra $\cC$ endowed with the action of $\Delta$ the {\em Higgs-Tate algebra} associated with the deformation $\tX$. 

\subsection{}\label{tordef40}
The Higgs-Tate algebra $\cC$ is a model of Hyodo's ring (\cite{agt} II.15.6). 
We introduce a ``weak'' $p$-adic completion that serves as a period ring for the $p$-adic Simpson correspondence. 
For any rational number $r\geq 0$, we denote by $\cF^{(r)}$ the
$\hoR$-representation of $\Delta$ deduced from $\cF$ by inverse image
under the morphism of multiplication by $p^r$ on
$\txi^{-1}\Omega^1_{R/\co_K}\otimes_R\hoR$, so that we have an exact sequence
\begin{equation}\label{tordef40a}
0\rightarrow \hoR\rightarrow \cF^{(r)}\rightarrow \txi^{-1}\Omega^1_{R/\co_K}\otimes_R\hoR\rightarrow 0.
\end{equation}
We consider the $\hoR$-algebra
\begin{equation}\label{tordef40b}
\cC^{(r)}=\underset{\underset{n\geq 0}{\longrightarrow}}\lim\ \Sym^n_{\hoR}(\cF^{(r)}),
\end{equation}
where the transition morphisms are defined by mapping $x_1\otimes\dots \otimes x_n$ to $1\otimes x_1\otimes\dots \otimes x_n$. 
The action of $\Delta$ on $\cF^{(r)}$ induces an action on $\cC^{(r)}$ by ring
automorphisms, compatible with its action on $\hoR$.  
We denote by $\hcC^{(r)}$ its $p$-adic completion of $\cC^{(r)}$.

For all rational numbers $r'\geq r\geq 0$, we have an injective and
$\Delta$-equivariant canonical $\hoR$-homomorphism
$\alpha^{r,r'}\colon \cC^{(r')}\rightarrow \cC^{(r)}$. One easily
verifies that the induced homomorphism
$\halpha^{r,r'}\colon\hcC^{(r')}\rightarrow \hcC^{(r)}$ is
injective. We set
\begin{equation}\label{tordef40c}
\hcC^{(r+)}=\underset{\underset{t\in \mQ_{>r}}{\longrightarrow}}{\lim} \hcC^{(t)},
\end{equation}
which we identify with a sub-$\hoR$-algebra of $\hcC=\hcC^{(0)}$.  The
group $\Delta$ acts naturally on $\hcC^{(r+)}$ by ring automorphisms,
in a manner compatible with its actions on $\hoR$ and on $\hcC$.

We denote by  
\begin{equation}\label{tordef40d}
d_{\cC^{(r)}}\colon \cC^{(r)}\rightarrow \txi^{-1}\Omega^1_{R/\co_K}\otimes_R\cC^{(r)}
\end{equation}
the universal $\hoR$-derivation of $\cC^{(r)}$ (\cite{ag2} 3.2.19) and by
\begin{equation}\label{tordef40e}
d_{\hcC^{(r)}}\colon \hcC^{(r)}\rightarrow \txi^{-1}\Omega^1_{R/\co_K}\otimes_R\hcC^{(r)}
\end{equation}
its extension to $p$-adic completions. These derivations are clearly $\Delta$-equivariant, 
and are Higgs $\hoR$-fields with coefficients in $\txi^{-1}\Omega^1_{R/\co_K}$ (\cite{ag2} 3.2.19). 

For all rational numbers $r'\geq r\geq 0$, we have
\begin{equation}\label{tordef40f}
p^{r'}(\id \times \alpha^{r,r'}) \circ d_{\cC^{(r')}}=p^rd_{\cC^{(r)}}\circ \alpha^{r,r'}.
\end{equation}
Therefore, the derivations $(p^td_{\hcC^{(t)}})_{t\in \mQ_{>r}}$ induce an $\hoR$-derivation 
\begin{equation}\label{tordef40g}
d^{(r)}_{\hcC^{(r+)}}\colon \hcC^{(r+)}\rightarrow \txi^{-1}\Omega^1_{R/\co_K} \otimes_R\hcC^{(r+)},
\end{equation}
which is also the restriction of $p^rd_{\hcC^{(r)}}$ to $\hcC^{(r+)}$.

\subsection{}\label{tordef41}
For any $\hoR$-representation $M$ of $\Delta$, we denote by $\mH(M)$ the $\hRun$-module defined by
\begin{equation}\label{tordef41a}
\mH(M)=(M\otimes_{\hoR}\hcC^{(0+)})^\Delta.
\end{equation}
We equip it with the Higgs $\hRun$-field with coefficients in $\txi^{-1}\Omega^1_{R/\co_K}$ induced by $d^{(0)}_{\hcC^{(0+)}}$ \eqref{tordef40g}.   
We thus define a functor 
\begin{equation}\label{tordef41b}
\mH\colon \bRep_{\hoR}(\Delta) \rightarrow \bHM(\hRun,\txi^{-1}\Omega^1_{R/\co_K}).
\end{equation}

\subsection{}\label{tordef42}
For any Higgs $\hRun$-module $(N,\theta)$ with coefficients in $\txi^{-1}\Omega^1_{R/\co_K}$, we denote by $\mV(N)$ the $\hoR$-module defined by 
\begin{equation}\label{tordef42a}
\mV(N)=(N\otimes_{\hRun}\hcC^{(0+)})^{\theta_\tot=0},
\end{equation}
where $\theta_\tot=\theta\otimes \id+\id\otimes d^{(0)}_{\hcC^{(0+)}}$ is the total Higgs $\hRun$-field on $N\otimes_{\hRun}\hcC^{(0+)}$.
We equip it with the $\hoR$-semi-linear action of $\Delta$ induced by its natural action on $\hcC^{(0+)}$. 
We thus define a functor 
\begin{equation}\label{tordef42b}
\mV\colon \bHM(\hRun,\txi^{-1}\Omega^1_{R/\co_K})\rightarrow \bRep_{\hoR}(\Delta).
\end{equation}

\begin{defi}[\cite{ag2} 3.3.9] \label{tordef43}
We say that an $\hoR[\frac 1 p]$-representation $M$ of $\Delta$ is {\em Dolbeault} 
if the following conditions are satisfied:
\begin{itemize}
\item[(i)] $\mH(M)$ is a projective $\hRun[\frac 1 p]$-module of finite type;
\item[(ii)] the canonical morphism 
\begin{equation}\label{tordef43a}
\mH(M) \otimes_{\hRun}\hcC^{(0+)}\rightarrow  M\otimes_{\hoR}\hcC^{(0+)}
\end{equation}
is an isomorphism.
\end{itemize}
\end{defi}

We prove that any Dolbeault $\hoR[\frac 1 p]$-representation de $\Delta$ is continuous for the $p$-adic topology (\cite{ag2} 3.3.12). 

\begin{defi}[\cite{ag2} 3.3.10] \label{tordef44}
We say that a Higgs $\hRun[\frac 1 p]$-module $(N,\theta)$ with coefficients in $\txi^{-1}\Omega^1_{R/\co_K}$
is {\em solvable} if the following conditions are satisfied:
\begin{itemize}
\item[(i)] $N$ is a projective $\hRun[\frac 1 p]$-module of finite type type;
\item[(ii)] the canonical morphism
\begin{equation}\label{tordef44a}
\mV(N) \otimes_{\hoR}\hcC^{(0+)}\rightarrow  N\otimes_{\hRun}\hcC^{(0+)}
\end{equation}
is an isomorphism.
\end{itemize}
\end{defi}

The notions \ref{tordef43} and \ref{tordef44} 
do not depend on the choice of the $\tS$-deformation $\tX$ (\cite{ag2} 3.3.8). But they depend a priori on the absolute or relative case considered in \ref{tordef21}.
They are equivalent to the {\em smallness} conditions introduced by Faltings \cite{faltings3}, {\em i.e.}, triviality conditions modulo prescribed powers of $p$ (\cite{ag2} 3.4.29--3.4.31). 

\begin{prop}[\cite{ag2} 3.3.16] \label{tordef45}
The functors $\mH$ \eqref{tordef41a} and $\mV$ \eqref{tordef42a} induce equivalences of categories quasi-inverse to each other,
between the category of  Dolbeault $\hoR[\frac 1 p]$-representations of $\Delta$ and that of 
solvable Higgs $\hRun[\frac 1 p]$-modules with coefficients in $\txi^{-1}\Omega^1_{R/\co_K}$.
\end{prop}

\begin{prop}[\cite{ag2} 3.3.17]\label{htls30}
For any Dolbeault $\hoR[\frac 1 p]$-representation $M$ of $\Delta$,  
there exists a canonical functorial isomorphism in $\bD^+(\bMod(\hRun[\frac 1 p]))$ 
\begin{equation}
\rC_\cont^\bullet(\Delta, M)\stackrel{\sim}{\rightarrow} \mK^\bullet(\mH(M)),
\end{equation}
where $\rC_\cont^\bullet(\Delta, M)$ is the continuous cochain complex of $\Delta$ with values in $M$
and $\mK^\bullet(\mH(M))$ is the Dolbeault complex of the Higgs $\hRun[\frac 1 p]$-module $\mH(M)$ associated to $M$ \eqref{tordef41b}. 
\end{prop}

\begin{prop}[\cite{ag2} 3.5.5]\label{htls3}
Let $M$ be a projective $\hoR[\frac 1 p]$-module of finite type, equipped with an $\hoR[\frac 1 p]$-semi-linear action of $\Delta$. 
Then, the following properties are equivalent:
\begin{itemize}
\item[{\rm (i)}] The $\hoR[\frac 1 p]$-representation $M$ of $\Delta$ is Dolbeault 
and the associated Higgs $\hRun[\frac 1 p]$-module $(\mH(M),\theta)$ \eqref{tordef41a} is nilpotent (i.e., there exists a finite decreasing filtration 
$(\mH_i)_{0\leq i\leq n}$ of $\mH(M)$ by sub-$\hRun[\frac 1 p]$-modules such that $\mH_0=\mH(M)$, $\mH_n=0$ and for any $0\leq i\leq n-1$, we have
\begin{equation}
\theta(\mH_i)\subset \txi^{-1}\Omega^1_{R/\co_K}\otimes_R\mH_{i+1}).
\end{equation}
\item[{\rm (ii)}] There exists a projective $\hRun[\frac 1 p]$-module of finite type $N$, a Higgs $\hRun[\frac 1 p]$-field $\theta$ on $N$ with coefficients 
in $\txi^{-1}\Omega^1_{R/\co_K}$ and a $\cC$-linear and $\Delta$-equivariant isomorphism
of  Higgs $\hoR[\frac 1 p]$-modules 
\begin{equation}
N \otimes_{\hRun}\cC\stackrel{\sim}{\rightarrow}  M\otimes_{\hoR}\cC.
\end{equation}
\end{itemize}
Moreover, under these conditions, we have an isomorphism of Higgs $\hRun[\frac 1 p]$-modules 
\begin{equation}
\mH(M)\stackrel{\sim}{\rightarrow}  (N,\theta).
\end{equation}
\end{prop}

\begin{defi}[\cite{ag2} 3.5.6]\label{htls4}
We say that an $\hoR[\frac 1 p]$-representation $M$ of $\Delta$ is {\em Hodge-Tate} if it satisfies the equivalent conditions of \ref{htls3}.
\end{defi}

This notion does not depend on the choice of the $\tS$-deformation $\tX$, not even on the absolute or relative case \ref{tordef21}.

\begin{rema}
Tsuji \cite{tsuji5} developed an arithmetic version of the local $p$-adic Simpson correspondence.
He associates to a $p$-adic $\hoR[\frac 1 p]$-representation of $\Gamma=\pi_1(X^\star_\eta,\oy)$ \eqref{tordef3} 
a Higgs field and an arithmetic Sen operator satisfying a compatibility relation that forces the Higgs field to be nilpotent (see also \cite{he2}). 
This explains the relation with the work of Liu and Zhu \cite{lz}. The reader should beware, however, that our notion of Hodge-Tate local systems 
does not correspond to that of Liu and Zhu since we consider geometric local systems while they consider arithmetic local systems.
\end{rema}

\section{Global theory. Dolbeault modules}\label{dolb}

\subsection{}
The local $p$-adic Simpson correspondence described in §~\ref{tordef} cannot be easily glued into a global correspondence for schemes which are not small affine. 
Instead, we sheafify the Higgs-Tate algebra in Faltings topos and use it to build a global $p$-adic Simpson correspondence parallel to the local picture. 
Then, we prove by cohomological descent that the global correspondence is equivalent to the local one for small affine schemes. 

\subsection{}\label{dolb1}
Let $X$ be a smooth $S$-scheme. 
With the notation of \ref{tordef21}, we assume that there exists a smooth $\tS$-deformation $\tX$ of $X\otimes_{\co_K}\co_C$ that we fix in this section: 
\begin{equation}\label{dolb1a}
\xymatrix{
{X\otimes_{\co_K}\co_C}\ar[r]\ar[d]\ar@{}[rd]|\Box&{\tX}\ar[d]\\
{\Spec(\co_C)}\ar[r]&{\tS}}
\end{equation}
Note that the condition is superfluous in the {\em relative case}; we can indeed take $\tX=X\times_S\tS$.

We take again the notation of section \ref{ft}. Recall in particular that we have the ring $\ocB=\{U\mapsto \ocB_U\}$ of Faltings topos $\tE$ \eqref{ft5}. For every $n\geq 0$, 
we set $\ocB_n=\ocB/p^n\ocB$ and for every $U\in \ob(\Et_{/X})$, $\ocB_{U,n}=\ocB_U/p^n\ocB_U$, which is a ring of $U_{\oeta,\fet}$. 

For any small affine étale $X$-scheme $U$, there exists a canonical exact sequence of $\ocB_{U,n}$-modules of $U_{\oeta,\fet}$
\begin{equation}\label{dolb1b}
0\rightarrow \ocB_{U,n}\rightarrow \cF_{U,n}\rightarrow \txi^{-1}\Omega^1_{X/S}(U)\otimes_{\co_X(U)}\ocB_{U,n} \rightarrow 0, 
\end{equation}
such that for any geometric point $\oy$ of $U_\oeta$, we have a canonical isomorphism of $\ocB_{U,\oy}$-representations of $\pi_1(U_\oeta,\oy)$
\begin{equation}
(\cF_{U,n})_\oy\stackrel{\sim}{\rightarrow}\cF^\oy_U/p^n\cF^\oy_U,
\end{equation}
where $\cF^\oy_U$ is the Higgs-Tate $\ocB_{U,\oy}$-extension \eqref{tordef4c} defined relatively to the restriction of $\tX$ over $U$ (\cite{ag2} 4.4.5). 
For every rational number $r\geq 0$, let 
\begin{equation}\label{dolb1c}
0\rightarrow \ocB_{U,n}\rightarrow \cF_{U,n}^{(r)}\rightarrow \txi^{-1}\Omega^1_{X/S}(U)\otimes_{\co_X(U)}\ocB_{U,n} \rightarrow 0
\end{equation}
be the extension of $\ocB_{U,n}$-modules of $U_{\oeta,\fet}$ obtained from $\cF_{U,n}$ by pull-back
by the multiplication by $p^r$ on $\txi^{-1}\Omega^1_{X/S}(U)\otimes_{\co_X(U)}\ocB_{U,n}$, and let  
\begin{equation}\label{dolb1d}
 \cC^{(r)}_{U,n}=\underset{\underset{m\geq 0}{\longrightarrow}}\lim\ \Sym^m_{\ocB_{U,n}}(\cF^{(r)}_{U,n})
\end{equation}
be the associated $\ocB_{U,n}$-algebra of $U_{\oeta,\fet}$, where the transition morphisms are locally 
defined by mapping $x_1\otimes\dots \otimes x_m$ to $1\otimes x_1\otimes\dots \otimes x_m$.  

The formation of $\cF^{(r)}_{U,n}$ being functorial in $U$, the correspondences 
\begin{equation}\label{dolb1e}
\{U\mapsto \cF^{(r)}_{U,n}\} \ \ \ {\rm and}\ \ \ \{U\mapsto \cC^{(r)}_{U,n}\}
\end{equation}
define presheaves on the subcategory $E^\sm$ of $E$ of objects $(V\rightarrow U)$ such that $U$ is small affine \eqref{ft1}. 
The latter is topologically generating of $E$. 
Therefore, by taking associated sheaves, we get a $\ocB_n$-module $\cF^{(r)}_n$ and a $\ocB_n$-algebra $\cC^{(r)}_n$ of $\tE$ (\cite{ag2} 4.4.10).

\subsection{}\label{dolb3}
We set $\cS=\Spf(\co_C)$ and denote by $\fX$ the formal $p$-adic completion of $\oX=X\times_S\oS$. 
Similarly, to take into account the $p$-adic topology, we consider the formal $p$-adic completion of the ringed topos $(\tE,\ocB)$. 
First, we define the {\em special fiber} $\tE_s$ of $\tE$, a topos that fits into a commutative diagram
\begin{equation}\label{dolb3a}
\xymatrix{
{\tE_s}\ar[r]^{\sigma_s}\ar[d]_{\delta}&{X_{s,\et}}\ar[d]^{\iota}\\
{\tE}\ar[r]^\sigma&{X_\et}}
\end{equation}
where $\iota$ is the canonical injection (\cite{ag2} 4.3.3). 
Concretely, $\tE_s$ is the full subcategory of $\tE$ of objects $F$ such that $F|\sigma^*(X_\eta)$ 
is the final object of $\tE_{/\sigma^*(X_\eta)}$, and $\delta_*\colon \tE_s\rightarrow \tE$ is the canonical injection functor.

For every integer $n\geq 0$, $\ocB_n$ is an object of $\tE_s$. We denote by $\oX_n$ and $\oS_n$ the reductions of $\oX$ and $\oS$ modulo $p^n$. 
Then, the morphism $\sigma_s$ is underlying a canonical morphism of ringed topos
\begin{equation}\label{dolb3b}
\sigma_n\colon (\tE_s,\ocB_n)\rightarrow (X_{s,\et},\co_{\oX_n}),
\end{equation}
where we identified the étale topos of $X_s$ and $\oX_n$, since $k$ is algebraically closed.   

The formal $p$-adic completion of $(\tE,\ocB)$ is the ringed topos $(\tE_s^{\mN^\circ},\bvocB)$, where $\tE_s^{\mN^\circ}$ is the topos of projective systems of objects of $\tE_s$ 
indexed by the ordered set $\mN$ (\cite{agt} III.7) and $\bvocB=(\ocB_{n})_{n\geq 0}$. The morphisms $\sigma_n$ induce a morphism of topos 
\begin{equation}\label{dolb3e}
\hupsigma \colon (\tE_s^{\mN^\circ},\bvocB)\rightarrow (X_{s,\zar},\co_\fX).
\end{equation}

We work in the category $\bMod_{\mQ}(\bvocB)$ of {\em $\bvocB$-modules up to isogeny}, {\em i.e.}, the category having for objects $\bvocB$-modules,
and for any $\bvocB$-modules $\cF$ and $\cG$, 
\begin{equation}\label{dolb3f}
\Hom_{\bMod_{\mQ}(\bvocB)}(\cF,\cG)=\Hom_{\bMod(\bvocB)}(\cF,\cG)\otimes_\mZ\mQ.
\end{equation} 
We denote the localization functor $\bMod(\bvocB)\rightarrow \bMod_{\mQ}(\bvocB)$ by $\cF\mapsto \cF_\mQ$. 
We call {\em $\bvocB_\mQ$-modules} the objects of $\bMod_{\mQ}(\bvocB)$. 

\subsection{}
For any rational number $r\geq 0$ and any integer $n\geq 0$, we have a canonical locally split exact sequence
\begin{equation}\label{dolb3c}
0\rightarrow \ocB_n\rightarrow \cF^{(r)}_n\rightarrow 
\sigma^*_n(\txi^{-1}\Omega^1_{\oX_n/\oS_n})\rightarrow 0
\end{equation}
and a canonical isomorphism of $\ocB_n$-algebras  
\begin{equation}\label{dolb3d}
\cC^{(r)}_n \stackrel{\sim}{\rightarrow}\underset{\underset{m\geq 0}{\longrightarrow}}\lim\ \Sym^m_{\ocB_n}(\cF^{(r)}_n),
\end{equation}
where the transition morphisms are locally defined by mapping $x_1\otimes\dots \otimes x_m$ to $1\otimes x_1\otimes\dots \otimes x_m$ (\cite{ag2} 4.4.11).
We set $\bvcF^{(r)}=(\cF^{(r)}_n)_{n\geq 0}$ which is a $\bvocB$-module and $\bvcC^{(r)}=(\cC^{(r)}_n)_{n\geq 0}$ which is a $\bvocB$-algebra.
We have a canonical exact sequence of $\bvocB$-modules 
\begin{equation}\label{dolb3g}
0\rightarrow \bvocB\rightarrow \bvcF^{(r)}\rightarrow 
\hupsigma^*(\txi^{-1}\Omega^1_{\fX/\cS})\rightarrow 0.
\end{equation}
The universal $\bvocB$-derivation of $\bvcC^{(r)}$ can be identified with a derivation
\begin{equation}\label{dolb3h}
d_{\bvcC^{(r)}}\colon \bvcC^{(r)}\rightarrow \hupsigma^*(\txi^{-1}\Omega^1_{\fX/\cS})\otimes_{\bvocB}\bvcC^{(r)}.
\end{equation}
It is a Higgs $\bvocB$-field. We denote by $\mK^\bullet(\bvcC^{(r)})$ the Dolbeault complex of the Higgs $\bvocB$-module $(\bvcC^{(r)},p^rd_{\bvcC^{(r)}})$.

For any rational numbers $r\geq r'\geq 0$, we have a canonical homomorphism of $\bvocB$-algebras $\bvcC^{(r)}\rightarrow \bvcC^{(r')}$. Observe that the restriction of 
the derivation $p^{r'}d_{\bvcC^{(r')}}$ is $p^{r}d_{\bvcC^{(r)}}$. Hence, we have a morphism of complexes 
\begin{equation}\label{dolb3i}
\mK^\bullet(\bvcC^{(r)})\rightarrow \mK^\bullet(\bvcC^{(r')}).
\end{equation}

\begin{prop}[\cite{agt} III.11.18, \cite{ag2} 4.4.32]\label{dolb4}
The canonical homomorphism
\begin{equation}\label{dolb4a}
\co_{\fX}[\frac 1 p]\rightarrow \underset{\underset{r\in \mQ_{>0}}{\longrightarrow}}{\lim}\  \hupsigma_*(\bvcC^{(r)})[\frac 1 p]
\end{equation}
is an isomorphism, and for any $q\geq 1$,
\begin{equation}\label{dolb4b}
\underset{\underset{r\in \mQ_{>0}}{\longrightarrow}}{\lim}\ \rR^q\hupsigma_*(\bvcC^{(r)})[\frac 1 p] =0.
\end{equation}
\end{prop}

This result is a sheafification of the computation of the Galois cohomology of the Higgs-Tate algebra over a small affine scheme (\cite{agt} II.12.5). 
The Galois cohomology computation relies on Faltings' almost purity result and the sheafification requires to prove a version modulo $p^n$,
up to some bounded defect (\cite{agt} II.12.7, \cite{ag2} 3.2.24).

\begin{prop}[\cite{agt} III.11.24, \cite{ag2} 4.4.36]\label{dolb5}
The canonical morphism of $\bMod_\mQ(\bvocB)$
\begin{equation}\label{dolb5a}
\bvocB_\mQ\rightarrow \underset{\underset{r\in \mQ_{>0}}{\longrightarrow}}{\lim}\ 
\rH^0(\mK^\bullet_\mQ(\bvcC^{(r)}))
\end{equation}
is an isomorphism, and for any $q\geq 1$, 
\begin{equation}\label{dolb5b}
\underset{\underset{r\in \mQ_{>0}}{\longrightarrow}}{\lim}\ \rH^q(\mK^\bullet_\mQ(\bvcC^{(r)}))=0.
\end{equation}
\end{prop}

This  result is a sheafification of the computation of the de Rham cohomology of the Higgs-Tate algebra over a small affine scheme (\cite{agt} II.12.3, \cite{ag2} 3.2.22). 

\subsection{}\label{dolb50}
Filtered inductive limits are not a priori representable in $\bMod_\mQ(\bvocB)$. However, we can naturally embed this category into the abelian category 
$\bIndMod(\bvocB)$ of ind-$\bvocB$-modules where filtered inductive limits are representable and which has better properties (\cite{ks2}, \cite{ag2} §~2.6 and §~2.7). 
In the same way, we can naturally embed the category of coherent $\co_\fX[\frac 1 p]$-modules into the category $\bIndMod(\co_\fX)$ of ind-$\co_\fX$-modules (\cite{ag2} 4.3.16).
The morphism $\hupsigma$ \eqref{dolb3e}  induces two adjoint functors 
\begin{equation}\label{dolb5c}
\xymatrix{
{\bIndMod(\bvocB)}\ar@<1ex>[r]^-(0.5){\rI \hupsigma_*}&{\bIndMod(\co_\fX)}\ar@<1ex>[l]^-(0.5){\rI \hupsigma^*}}
\end{equation}
that extend the adjoint functors $\hupsigma^*$ and $\hupsigma_*$. 

\begin{defi}
We call {\em Higgs $\co_\fX[\frac 1 p]$-bundle with coefficients in $\txi^{-1}\Omega^1_{\fX/\cS}$} 
any locally projective $\co_\fX[\frac 1 p]$-module of finite type $\cN$ (\cite{ag2} 2.1.11)
equipped with an $\co_\fX$-linear morphism $\theta\colon \cN\rightarrow \txi^{-1}\Omega^1_{\fX/\cS}\otimes_{\co_\fX}\cN$ such that $\theta\wedge \theta=0$.
\end{defi}

Observe that since the stalks of the ring $\co_\fX[\frac 1 p]$ are not necessarily local rings, 
locally projective $\co_\fX[\frac 1 p]$-modules of finite type are not necessarily locally free. 

\begin{defi}[\cite{ag2} 4.5.4]\label{dolb6}
Let $\cM$ be an ind-$\bvocB$-module, $\cN$ a Higgs $\co_\fX[\frac 1 p]$-bundle with coefficients in $\txi^{-1}\Omega^1_{\fX/\cS}$. 
\begin{itemize}
\item[(i)] We say that $\cM$ and  $\cN$ are {\em $r$-associated} (for $r\in \mQ_{>0}$) if there exists an isomorphism of ind-$\bvcC^{(r)}$-modules 
\begin{equation}\label{dolb6a}
\cM\otimes_{\bvocB}\bvcC^{(r)}\stackrel{\sim}{\rightarrow}\rI\hupsigma^*(\cN)\otimes_{\bvocB}\bvcC^{(r)},
\end{equation}
compatible with the total Higgs $\bvocB$-fields with coefficients in $\hupsigma^*(\txi^{-1}\Omega^1_{\fX/\cS})$, 
where $\cM$ is equipped with the zero Higgs field and $\bvcC^{(r)}$ with $p^rd_{\bvcC^{(r)}}$. 
\item[(ii)] We say that $\cM$ and  $\cN$ are {\em associated} if they are $r$-associated for a rational number $r>0$.
\end{itemize}
\end{defi}

In fact, \eqref{dolb6a} is an isomorphism of {\em ind-$\bvcC^{(r)}$-modules with $p^r$-connection relatively to 
the extension $\bvcC^{(r)}/\bvocB$} (\cite{ag2} 4.5.1).
In particular, for any rational numbers $r\geq r'>0$, if $\cM$ and $\cN$ are $r$-associated, they are $r'$-associated.

\begin{defi}[\cite{ag2} 4.5.5]\label{dolb7}
\
\begin{itemize}
\item[(i)] We say that an ind-$\bvocB$-module is {\em Dolbeault} if it is associated to a  
Higgs $\co_\fX[\frac 1 p]$-bundle with coefficients in $\txi^{-1}\Omega^1_{\fX/\cS}$.
\item[(ii)] We say that a Higgs $\co_\fX[\frac 1 p]$-bundle with coefficients in $\txi^{-1}\Omega^1_{\fX/\cS}$ is {\em solvable}
if it is associated to an ind-$\bvocB$-module. 
\end{itemize}
\end{defi}

The property for an ind-$\bvocB$-module to be Dolbeault does not depend on the choice of the deformation $\tX$ \eqref{dolb1}
provided that we stay in one of the settings, absolute or relative \eqref{tordef21} (\cite{ag2} 4.10.4).
The property for a Higgs $\co_\fX[\frac 1 p]$-bundle to be solvable depends a priori on the deformation $\tX$ (see however \cite{ag2} 4.10.12).

The notion of being Dolbeault applies to $\bvocB_\mQ$-modules \eqref{dolb50}. We call the associated Higgs bundles {\em rationally solvable}.

In (\cite{agt} III.12.11), we considered only Dolbeault $\bvocB_\mQ$-modules and we requested moreover that they are adic of finite type. 
We renamed them in (\cite{ag2} 4.6.6) {\em strongly Dolbeault} $\bvocB_\mQ$-modules and renamed the associated Higgs bundles {\em strongly solvable}. 
The finiteness condition is important for matching the local and global theories for small affine schemes \eqref{htls8}.

\begin{teo}[\cite{ag2} 4.5.20]\label{dolb8}
There are explicit equivalences of categories quasi-inverse to each other 
\begin{equation}\label{dolb8a}
\xymatrix{
{\bIndMod^\Dolb(\bvocB)}\ar@<1ex>[r]^-(0.5){\cH}&{\bHM^\sol(\co_\fX[\frac 1 p], \txi^{-1}\Omega^1_{\fX/\cS})}
\ar@<1ex>[l]^-(0.5){\cV}}
\end{equation}
between the category of Dolbeault  ind-$\bvocB$-modules and the category of solvable Higgs $\co_\fX[\frac 1 p]$-bundles with coefficients in $\xi^{-1}\Omega^1_{\fX/\cS}$. 
\end{teo}

These functors are explicitly defined (\cite{ag2} 4.5.7). Let $\vupsigma_*$ be the composed functor
\begin{equation}\label{dolb8d}
\xymatrix{\bIndMod(\bvocB)\ar[r]^-(0.5){\rI \hupsigma_*}&{\bIndMod(\co_\fX)}\ar[r]^-(0.5){\kappa_{\co_\fX}}&{\bMod(\co_\fX)}},
\end{equation}
where $\rI \hupsigma_*$ is defined in \eqref{dolb5c} and 
\begin{equation}\label{dolb8c}
\kappa_{\co_\fX}(\indcolim\alpha)= \underset{\underset{J}{\longrightarrow}}{\lim}\ \alpha.
\end{equation}
Then, the functor $\cH$ can in fact be defined for any ind-$\bvocB$-module $\cM$ by
\begin{equation}\label{dolb8e}
\cH(\cM)=\underset{\underset{r\in \mQ_{>0}}{\longrightarrow}}{\lim}\ \vupsigma_*(\cM\otimes_{\bvocB}\bvcC^{(r)},p^r\id\otimes d_{\bvcC^{(r)}}).
\end{equation} 
We have a similar definition of $\cV$. The functors $\cH$ and $\cV$ depend a priori on the deformation $\tX$.

\begin{teo}[\cite{ag2} 4.7.4]\label{dolb9}
For any Dolbeault ind-$\bvocB$-module $\cM$ and any integer $q\geq 0$, there is a canonical functorial isomorphism of $\bD^+(\bMod(\co_\fX))$
\begin{equation}\label{dolb9a}
\rR\vupsigma_*(\cM)\stackrel{\sim}{\rightarrow}\mK^\bullet(\cH(\cM)),
\end{equation}
where $\mK^\bullet(\cH(\cM))$ is the Dolbeault complex of $\cH(\cM)$.
\end{teo}

This is the global analogue of \ref{htls30}.

\subsection{}\label{dolb90}
The morphism $\psi$ \eqref{ft3a} induces a morphism of topos
\begin{equation}
\bvpsi\colon X_{\oeta,\et}^{\mN^\circ}\rightarrow \tE^{\mN^\circ},
\end{equation}
where $X_{\oeta,\et}^{\mN^\circ}$ is the topos of projective systems of objects of $X_{\oeta,\et}$, indexed by the ordered set $\mN$. 
We denote by $\bvmZ_p$ the ring $(\mZ/p^n\mZ)_{n\geq 0}$ of $X_{\oeta,\et}^{\mN^\circ}$. 

We say that a $\bvmZ_p$-module $M=(M_n)_{n\in \mN}$ of $X_{\oeta,\et}^{\mN^\circ}$ is a {\em local system} if the following two conditions are satisfied:
\begin{itemize}
\item[(a)] $M$ is $p$-adic, {\em i.e.}, for any integers $n\geq m\geq 0$, the morphism $M_n/p^mM_n\rightarrow M_m$ deduced from the transition morphism 
$M_n\rightarrow M_m$, is an isomorphism; 
\item[(b)]  for any integer $n\geq 0$, the $\mZ/p^n\mZ$-module $M_n$ of $X_{\oeta,\et}$ is locally constant constructible. 
\end{itemize}

\begin{cor}[\cite{ag2} 4.12.6]\label{dolb10}
Let $M=(M_n)_{n\geq 0}$ be a $\bvmZ_p$-local system of $X_{\oeta,\et}^{\mN^\circ}$, $\cM=\bvpsi_*(M)\otimes_{\bvmZ_p}\bvocB$. 
Assume that $X$ is proper over $S$ and that the $\bvocB_\mQ$-module $\cM_\mQ$ is Dolbeault.
Then, there exists a canonical spectral sequence
\begin{equation}\label{dolb10a}
\rE_2^{i,j}=\rH^i(X_s,\rH^j(\mK^\bullet))\Rightarrow \rH^{i+j}(X_{\oeta,\et}^{\mN^\circ},M)\otimes_{\mZ_p}C,
\end{equation}
where $\mK^\bullet$ be the Dolbeault complex of $\cH(\cM_\mQ)$. 
\end{cor}

It follows from \ref{ft6} and \ref{dolb9}.

\begin{rema}
In \ref{dolb10}, if we take $M=\bvmZ_p$, then $\cM=\bvocB$, the $\bvocB_\mQ$-module $\bvocB_\mQ$ is Dolbeault and $\cH(\bvocB_\mQ)$ is equal to $\co_\fX[\frac 1 p]$ 
equipped with the zero Higgs field (\cite{ag2} 4.6.10). The spectral sequence \eqref{dolb10a} is the Hodge-Tate spectral sequence (\cite{ag1} 6.4.6). 
Observe that the construction \ref{dolb10} of this spectral sequence shows directly that it degenerates at $\rE_2$ and that
the abutment filtration is split without using Tate's theorem on the Galois cohomology of $C(j)$.
This construction applies in particular by taking for $\tX$ in the relative case \eqref{tordef21} the trivial deformation \eqref{dolb1}.
\end{rema}

\begin{defi}[\cite{ag2} 4.11.1]\label{htls5}
We call {\em Hodge-Tate $\bvocB_\mQ$-module} any Dolbeault $\bvocB_\mQ$-module $\cM$ \eqref{dolb7} 
whose associated Higgs $\co_\fX[\frac 1 p]$-bundle $(\cH(\cM),\theta)$ \eqref{dolb8a} is nilpotent, {\em i.e.}, there exists a finite 
decreasing filtration $(\cH_i(\cM))_{0\leq i\leq n}$ of $\cH(\cM)$ by coherent sub-$\co_\fX[\frac 1 p]$-modules 
such that $\cH_0(\cM)=\cH(\cM)$, $\cH_n(\cM)=0$ and that for any $0\leq i\leq n-1$, we have
\begin{equation}\label{htls5a}
\theta(\cH_i(\cM))\subset \txi^{-1}\Omega^1_{\fX/\cS}\otimes_{\co_\fX}\cH_{i+1}(\cM).
\end{equation}
\end{defi} 

This notion does not depend on the choice of the $\tS$-deformation $\tX$, not even on the absolute or relative case \eqref{tordef21} (\cite{ag2} 4.11.9). 

\begin{prop}[\cite{ag2} 4.11.2]\label{htls6}
The functors $\cH$ and $\cV$ \eqref{dolb8a} induce equivalences of categories quasi-inverse to each other
\begin{equation}\label{htls6a}
\xymatrix{
{\bMod^\HT_\mQ(\bvocB)}\ar@<1ex>[r]^-(0.5){\cH}&{\bHM^\qsolnilp(\co_\fX[\frac 1 p], \txi^{-1}\Omega^1_{\fX/\cS})}
\ar@<1ex>[l]^-(0.5){\cV}}
\end{equation} 
between the category of Hodge-Tate $\bvocB_\mQ$-modules and the category of rationally solvable and 
nilpotent Higgs $\co_\fX[\frac 1 p]$-bundles with coefficients in $\txi^{-1}\Omega^1_{\fX/\cS}$ \eqref{dolb7}. 
\end{prop}

\subsection{}\label{htls7}
We assume in the remaining part of this section that $X$ is small affine \eqref{tordef3}, that $X_s\not=\emptyset$, and for simplicity 
that $X_\oeta$ is connected. We fix  a geometric point $\oy$ of $X_\oeta$. We set $\Delta=\pi_1(X_\oeta,\oy)$, and we denote 
by $\bB_\Delta$ the classifying topos of $\Delta$ and by
\begin{equation}\label{htls7a}
\nu\colon X_{\oeta,\fet} \stackrel{\sim}{\rightarrow} \bB_\Delta
\end{equation}
the fiber functor of $X_{\oeta,\fet}$ at $\oy$ (\cite{agt}  (VI.9.8.4)). We denote by $\upbeta$ the composed functor
\begin{equation}\label{htls7b}
\upbeta \colon \tE\rightarrow \bB_\Delta, \ \ \ F\mapsto \nu\circ (\beta_*(F)), 
\end{equation}
where $\beta$ is the morphism of topos \eqref{ft3c}. 
We thus define a functor from the category of abelian sheaves of $\tE$ into the category of $\mZ[\Delta]$-modules. 
The latter being left exact, we denote by $\rR^q \upbeta$ $(q\geq 0)$ its right derived functors. 
For any abelian sheaf $F$ of $\tE$ and any integer $q\geq 0$, we have a canonical functorial isomorphism
\begin{equation}
\rR^q\upbeta(F)\stackrel{\sim}{\rightarrow} \nu\circ (\rR^q\beta_*(F)).
\end{equation}

For any abelian sheaf $F=(F_n)_{n\geq 0}$ of $\tE^{\mN^\circ}$, we set
\begin{equation}\label{htls7c}
\hupbeta(F)=\underset{\underset{n\geq 0}{\longleftarrow}}{\lim}\ \upbeta(F_n).
\end{equation}
We thus define a functor from the category of abelian sheaves of $\tE^{\mN^\circ}$ into the category of $\mZ[\Delta]$-modules.
The latter being left exact, we abusively denote by $\rR^q\hupbeta(F)$ ($q\geq 0$) its right derived functors. 
By (\cite{jannsen} 1.6), we have a canonical exact sequence 
\begin{equation}\label{htls7d}
0\rightarrow \rR^1 \underset{\underset{n\geq 0}{\longleftarrow}}{\lim}\ \rR^{q-1}\upbeta(F_n)\rightarrow
\rR^q\hupbeta(F)\rightarrow 
\underset{\underset{n\geq 0}{\longleftarrow}}{\lim}\ \rR^q\upbeta(F_n)\rightarrow 0,
\end{equation}
where we set $\rR^{-1}\upbeta(F_n)=0$ for all $n\geq 0$. 

We set $\oR=\nu(\ocB_X)$, which is nothing but the algebra defined in \eqref{ht5c} equipped with the canonical action of $\Delta$. 
For any integer $q\geq 0$, $\rR^q\hupbeta$ induces a functor that we also denote by
\begin{equation}\label{htls7e}
\rR^q\hupbeta\colon \bMod(\bvocB)\rightarrow \bRep_{\hoR}(\Delta). 
\end{equation}
The latter induces a functor that we also denote by 
\begin{equation}\label{htls7f}
\rR^q\hupbeta\colon \bMod_\mQ(\bvocB)\rightarrow \bRep_{\hoR[\frac 1 p]}(\Delta). 
\end{equation}

\begin{teo}[\cite{ag2} 4.8.31]\label{htls8}
We keep the assumptions and notation of \ref{htls7}. Then,
\begin{itemize}
\item[{\rm (i)}] The functor $\hupbeta$ \eqref{htls7f} induces an equivalence of categories
\begin{equation}\label{htls8a}
\bMod^\sDolb_\mQ(\bvocB)\stackrel{\sim}{\rightarrow} \bRep_{\hoR[\frac 1 p]}^{\Dolb}(\Delta),
\end{equation} 
between the category of strongly Dolbeault $\bvocB_\mQ$-modules \eqref{dolb7} and the category of Dolbeault 
$\hoR[\frac 1 p]$-representations of $\Delta$ \eqref{tordef43}.
\item[{\rm (ii)}] For any strongly Dolbeault $\bvocB_\mQ$-module $\cM$ and any integer $q\geq 1$, we have
\begin{equation}\label{htls8b}
\rR^q\hupbeta(\cM)=0.
\end{equation}
\end{itemize}
\end{teo}

This follows from a cohomological descent result for strongly Dolbeault $\bvocB_\mQ$-modules (\cite{ag2} 4.8.16) which {\em in fine} reduces to 
a cohomological descent result for the ring $\bvocB_\mQ$ (\cite{ag1} 4.6.30). 

We show that under the equivalence of categories \eqref{htls8a}, the functors $\cH$ \eqref{dolb8a} and $\mH$ \eqref{tordef41b}
(resp. $\cV$ \eqref{dolb8a} and $\mV$ \eqref{tordef42b}) correspond (\cite{ag2} 4.8.32).

\section{\texorpdfstring{Functoriality of the $p$-adic Simpson correspondence by proper direct image}
{Functoriality of the p-adic Simpson correspondence by proper direct image}}\label{fpscpdi}

\subsection{}\label{fpscpdi1} 
Let $g\colon X'\rightarrow X$ be a smooth morphism of smooth $S$-schemes. We equip with a prime $^\prime$ the objects associated to $X'/S$. 
By functoriality of the Faltings topos, $g$ induces a canonical morphism $\Theta$ between Faltings topos that fits into a commutative diagram 
\begin{equation}\label{fpscpdi1a}
\xymatrix{
{X'_{\oeta,\et}}\ar[d]_{g_\oeta}\ar[r]^{\psi'}&{\tE'}\ar[d]^{\Theta}\ar[r]^{\sigma'}&{X'_\et}\ar[d]^g\\
{X_{\oeta,\et}}\ar[r]^\psi&{\tE}\ar[r]^{\sigma}&{X_\et}}
\end{equation}
We have also a canonical ring homomorphism 
\begin{equation}\label{fpscpdi1b}
\ocB\rightarrow \Theta_*(\ocB').
\end{equation}

\begin{teo}[\cite{faltings2}, \cite{ag1} 5.7.3] \label{fpscpdi2} 
Assume that $g\colon X'\rightarrow X$ is projective, and let  $F'$ be a locally constant constructible sheaf of 
$(\mZ/p^n\mZ)$-modules of $X'_{\oeta,\et}$ $(n\geq 1)$. Then, for any integer $i\geq 0$, the canonical morphism 
\begin{equation}\label{fpscpdi2a} 
\psi_*(\rR^ig_{\oeta*}(F'))\otimes_{\mZ_p}\ocB\rightarrow \rR^i\Theta_*(\psi'_*(F')\otimes_{\mZ_p}\ocB')
\end{equation}
is an almost isomorphism. 
\end{teo}

Observe that the sheaves $\rR^ig_{\oeta*}(F)$ are locally constant constructible on $X_\oeta$ 
by the smooth and the proper base change theorems. 

Faltings formulated this {\em relative version} of his main $p$-adic comparison theorem 
in \cite{faltings2} and he very roughly sketched a proof in the appendix. 
Some arguments have to be modified and our actual proof in \cite{ag1} requires much more work.

The condition {\em $g$ projective} could be replaced by the condition {\em $g$ proper} by using the cohomological descent established in \cite{he1}.

\subsection{}\label{fpscpdi3}
With the notation of \ref{tordef21}, we assume that there exists a commutative diagram with Cartesian squares
\begin{equation}\label{fpscpdi3b}
\xymatrix{
{X'\otimes_{\co_K}\co_C}\ar[d]_{g\otimes\id}\ar@{}[rd]|\Box\ar[r]&{\tX'}\ar[d]^{\tg}\\
{X\otimes_{\co_K}\co_C}\ar[r]\ar[d]\ar@{}[rd]|\Box&{\tX}\ar[d]\\
{\Spec(\co_C)}\ar[r]&{\tS}}
\end{equation}
where $\tX$ and $\tX'$ are smooth $\tS$-schemes, {\em that we fix in this section}. 
Note that the condition is superfluous in the {\em relative case}; we can indeed take $\tg=g\times_S\tS$.
Observe that $\tg$ is smooth by (\cite{ega4} 17.11.1).

Diagram \eqref{fpscpdi1a} induces a commutative diagram of morphisms of ringed topos
\begin{equation}\label{fpscpdi3a}
\xymatrix{
{(\tE'^{\mN^\circ}_s,\bvocB')}\ar[d]_{\hupsigma'}\ar[r]^{\bvuptheta}&{(\tE^{\mN^\circ}_s,\bvocB)}\ar[d]^-(0.5){\hupsigma}\\
{(X'_{s,\zar},\co_{\fX'})}\ar[r]^{\fgg}&{(X_{s,\zar},\co_{\fX})}}
\end{equation}
where the horizontal arrows are induced by $\Theta$ and $g$, and the vertical arrows are induced by $\sigma'$ and $\sigma$ \eqref{dolb3e}. 
We prove in (\cite{ag2} 6.3.20) that for any rational number $r\geq 0$, the lifting $\tg$ \eqref{fpscpdi3b} induces a homomorphism of the Higgs-Tate algebras 
\begin{equation}\label{fpscpdi3c}
\bvuptheta^*(\bvcC^{(r)})\rightarrow \bvcC'^{(r)},
\end{equation}
whose construction is rather subtle.

\begin{teo}[\cite{ag2} 6.5.24]\label{fpscpdi4}
Assume $g\colon X'\rightarrow X$ proper. 
Let $\cM$ be a Dolbeault ind-$\bvocB'$-module \eqref{dolb7},  
\begin{equation}\label{fpscpdi4a}
\cH'(\cM)\rightarrow \txi^{-1}\Omega^1_{\fX'/\cS}\otimes_{\co_{\fX'}}\cH'(\cM)
\end{equation}
the associated Higgs bundle \eqref{dolb8a}, 
\begin{equation}\label{fpscpdi4b}
\ucH'(\cM)\rightarrow \txi^{-1}\Omega^1_{\fX'/\fX}\otimes_{\co_{\fX'}}\ucH'(\cM)
\end{equation}
the relative Higgs bundle induced by \eqref{fpscpdi4a}, $\umK^\bullet$ the Dolbeault complex of $\ucH'(\cM)$. 
Then, for any integer $q\geq 0$, there exists a rational number $r>0$ and a $\bvcC^{(r)}$-isomorphism 
\begin{equation}\label{fpscpdi4c}
\rR^q\rI\bvuptheta_*(\cM)\otimes_{\bvocB}\bvcC^{(r)}\stackrel{\sim}{\rightarrow}
\rI\hupsigma^*(\rR^q\fgg_*(\umK^\bullet))\otimes_{\bvocB}\bvcC^{(r)},
\end{equation}
compatible with the total Higgs fields, where $\bvcC^{(r)}$ is equipped with the Higgs field $p^rd_{\bvcC^{(r)}}$, $\rR^q\rI\bvuptheta_*(\cM)$ 
with the zero Higgs field and $\rR^q\fgg_*(\umK^\bullet)$ with the Katz-Oda field {\rm (\cite{ag2} 2.5.16)}.
\end{teo}

Observe that the $\co_\fX[\frac 1 p]$-module $\rR^q\fgg_*(\umK^\bullet)$ is coherent, and that the functor 
\begin{equation}\label{fpscpdi4d}
\rI\hupsigma^*\colon \bMod^\coh(\co_{\fX}[\frac 1 p])\rightarrow \bIndMod(\bvocB)
\end{equation}
is exact (\cite{ag2} 6.2.6).

\begin{cor}[\cite{ag2} 6.5.25]\label{fpscpdi5}
Under the assumptions of \ref{fpscpdi4}, if the $\co_{\fX}[\frac 1 p]$-module $\rR^q\fgg_*(\umK^\bullet)$ is locally projective of finite type, then 
the ind-$\bvocB$-module $\rR^q\rI\bvuptheta_*(\cM)$ is {\em Dolbeault}, and we have an isomorphism
\begin{equation}\label{fpscpdi5a}
\cH(\rR^q\rI\bvuptheta_*(\cM))\stackrel{\sim}{\rightarrow} \rR^q\fgg_*(\umK^\bullet),
\end{equation}
where $\rR^q\fgg_*(\umK^\bullet)$ is equipped with the Katz-Oda field.
\end{cor}

\begin{cor}[\cite{ag2} 6.5.34]\label{fpscpdi8}
Assume $g\colon X'\rightarrow X$ projective, and let $\cM^n=\bvpsi_*(\rR^n\bvg_{\oeta*}(\bvmZ_p))\otimes_{\bvmZ_p}\bvocB$ 
for an integer $n\geq 0$.  
Then, the $\bvocB_\mQ$-module $\cM^n_\mQ$ is Hodge-Tate \eqref{htls5} and we have an isomorphism
\begin{equation}\label{fpscpdi8a}
\cH(\cM^n_\mQ)\stackrel{\sim}{\rightarrow} \oplus_{0\leq i\leq n}\rR^ig_*(\txi^{i-n}\Omega^{n-i}_{X'/X}) \otimes_{\co_X}\co_\fX[\frac 1 p],
\end{equation}
where the Higgs field on the right hand side is induced by the Kodaira-Spencer maps of $g$ 
\begin{equation}\label{fpscpdi8b}
\rR g_*(\txi^{-j}\Omega^j_{X'/X})\rightarrow \txi^{-1}\Omega^1_{X/S}\otimes_{\co_X}\rR g_*(\txi^{1-j}\Omega^{j-1}_{X'/X})[+1].
\end{equation}
\end{cor}

It follows from \ref{fpscpdi2} and \ref{fpscpdi5}. Indeed, $\bvpsi'_*(\bvmZ_p)=\bvmZ_p$, the $\bvocB'_\mQ$-module $\bvocB'_\mQ$ is Dolbeault and 
$\ucH'(\bvocB'_\mQ)$ is equal to $\co_{\fX'}[\frac 1 p]$ equipped with the zero Higgs field (\cite{agt} III.12.14). 
Therefore, with the notation of \ref{fpscpdi4}, for any $q\geq 0$, the $\co_{\fX}[\frac 1 p]$-module 
$\rR^q\fgg_*(\umK^\bullet)$ is locally free of finite type by (\cite{deligne1} 5.5), which completes the proof of the first statement. 
The second statement follows easily from the definition of the Katz-Oda field (\cite{katz2} 1.2, \cite{ag2} 2.5.17).

\begin{cor}[\cite{ag2} 6.5.28]\label{fpscpdi6}
Let $M=(M_n)_{n\geq 0}$ be a $\bvmZ_p$-local system of $X'^{\mN^\circ}_{\oeta,\et}$ \eqref{dolb90}.  
We set $\cM=\bvpsi'_*(M)\otimes_{\bvmZ_p}\bvocB'$, 
and we assume that $g\colon X'\rightarrow X$ is projective and that the $\bvocB'_\mQ$-module $\cM_\mQ$ is Dolbeault.
Then, there exists a rational number $r>0$ and a spectral sequence
\begin{equation}\label{fpscpdi6a}
\rE_2^{i,j}=\hupsigma^*_\mQ(\rR^i\fgg_*(\rH^j(\umK^\bullet)))\otimes_{\bvocB_\mQ}\bvcC^{(r)}_\mQ
\Rightarrow \bvpsi_*(\rR^{i+j}\bvg_{\oeta*}(M))\otimes_{\bvmZ_p}\bvcC^{(r)}_\mQ,
\end{equation}
where $\umK^\bullet$ is the Dolbeault complex of the relative Higgs $\co_{\fX'}[\frac 1 p]$-bundle $\ucH'(\cM_\mQ)$ \eqref{fpscpdi4b}. 
\end{cor}
 
It follows from \ref{fpscpdi2} and \ref{fpscpdi4}.

\begin{cor}[\cite{ag2} 6.5.29]\label{fpscpdi60}
Assume $g\colon X'\rightarrow X$ projective.
Then, there exists rational number $r>0$ and for any integer $n\geq 0$, a canonical isomorphism of $\bvcC^{(r)}_\mQ$-modules 
\begin{equation}\label{fpscpdi60a}
\bvpsi_*(\rR^n\bvg_{\oeta*}(\bvmZ_p))\otimes_{\bvmZ_p}\bvcC^{(r)}_\mQ
\stackrel{\sim}{\rightarrow}\oplus_{0\leq i\leq n}\sigma^*(\rR^ig_*(\Omega^{n-i}_{X'/X}))\otimes_{\sigma^*(\co_X)}\bvcC^{(r)}_\mQ(i-n).
\end{equation}
\end{cor}

It follows from \ref{fpscpdi4} applied to $\cM=\bvocB'_\mQ$, since the $\bvocB'_\mQ$-module $\bvocB'_\mQ$ is Dolbeault and 
$\ucH'(\bvocB'_\mQ)$ is the trivial bundle $\co_{\fX'}[\frac 1 p]$ equipped with the null Higgs field.  

This corollary applies in particular by taking for $\tg$ in the relative case \eqref{tordef21} the trivial deformation \eqref{fpscpdi3}.

\begin{teo}[\cite{ag1} 6.7.5]\label{fpscpdi7}
Assume $g\colon X'\rightarrow X$ projective. Then, we have a canonical spectral sequence of $\bvocB_\mQ$-modules, {\em the relative Hodge-Tate spectral sequence},
\begin{equation}\label{fpscpdi7a}
\rE_2^{i,j}=\sigma^*(\rR^ig_*(\Omega^j_{X'/X}))\otimes_{\sigma^*(\co_X)}\bvocB_\mQ(-j)\Rightarrow \bvpsi_*(\rR^{i+j}\bvg_{\oeta*}(\bvmZ_p))\otimes_{\bvmZ_p}\bvocB_\mQ.
\end{equation}
\end{teo}

This spectral sequence does not require any deformation \eqref{fpscpdi3b}. 
It is $G_K$-equivariant for the natural $G_K$-equivariant structures on the various topos and objects involved. 
Hence it degenerates at $\rE_2$. But the abutment filtration does not split in general. However, we can check that it splits 
after base change from $\bvocB$ to $\bvcC^{(r)}$ for a rational number $r>0$, 
and that it corresponds to the decomposition \eqref{fpscpdi60a}.

\subsection{}\label{fpscpdi9} 
The proof of \ref{fpscpdi4} can be divided into three steps. First, we compute the relative Galois and Higgs cohomologies of the Higgs-Tate algebra by adapting Faltings' computation 
in the absolute case. Second, to sheafify these computations, we consider the following fiber product of topos 
\begin{equation}\label{fpscpdi9a} 
\xymatrix{
\tE'\ar[rd]^{\sigma'}\ar[d]_{\tau}&\\
{\tE\times_{X_\et}X'_\et}\ar[r]\ar[d]\ar@{}|\Box[rd]&{X'_\et}\ar[d]^g\\
{\tE}\ar[r]^-(0.5){\sigma}&{X_\et}}
\end{equation}
The local relative Galois cohomology computation can be globalized into a computation of the sheaves $\rR^i\tau_*(\cC'^{(r)}_n)$. 
The last step is a base change statement relatively to the Cartesian square. 

It turns out that there is a very natural site underlying the topos $\tE\times_{X_\et}X'_\et$, which is a relative version of Faltings topos
and whose definition was inspired by oriented products of topos (beyond the covanishing topos which inspired the usual Faltings topos). 

\section{Relative Faltings topos}\label{rft} 

\subsection{}\label{rft1} 
Let $g\colon X'\rightarrow X$ be a morphism of $S$-schemes. 
We denote by $G$ the category of morphisms $(W\rightarrow U\leftarrow V)$
above the canonical morphisms $X'\rightarrow X\leftarrow X_\oeta$, {\em i.e.}, commutative diagrams 
\begin{equation}\label{rft1a}
\xymatrix{W\ar[r]\ar[d]&U\ar[d]&V\ar[l]\ar[d]\\
X'\ar[r]&X&X_\oeta\ar[l]}
\end{equation}
such that $W$ is étale over $X'$, $U$ is étale over $X$ and the canonical morphism $V\rightarrow U_\oeta$ is {\em finite étale}. 
We equip it with the topology generated by coverings 
\[
\{(W_i\rightarrow U_i\leftarrow V_i)\rightarrow (W\rightarrow U \leftarrow V)\}_{i\in I}
\]
of the following three types~:
\begin{itemize}
\item[(a)] $U_i=U$, $V_i=V$ for all $i\in I$ and $(W_i\rightarrow W)_{i\in I}$ is a covering;
\item[(b)] $W_i=W$, $U_i=U$ for all $i\in I$ and $(V_i\rightarrow V)_{i\in I}$ is a covering;
\item[(c)]  diagrams
\begin{equation}\label{rft1b}
\xymatrix{
W'\ar[r]\ar@{=}[d]&U'\ar[d]\ar@{}|\Box[rd]&V'\ar[l]\ar[d]\\
W\ar[r]&U&V\ar[l]}
\end{equation}
where $U'\rightarrow U$ is any morphism and the right square is Cartesian. 
\end{itemize}

We denote by $\tG$ the topos of sheaves of sets on $G$ and we call it the {\em relative Faltings topos} associated with the pair of morphisms 
$(X_\oeta\rightarrow X,X'\rightarrow X)$.  

There is a canonical morphism of topos 
\begin{equation}\label{rft1c}
\pi\colon \tG\rightarrow X'_\et, \ \ \ W\in \ob(\Et_{/X'}) \mapsto \pi^*(W)=(W\rightarrow X\leftarrow X_\oeta)^a.
\end{equation}

\subsection{}\label{rft2} 
If $X'=X$, the topos $\tG$ is canonically equivalent to the Faltings topos $\tE$. 
By functoriality of the relative Faltings topos, we then get a natural factorization of the canonical morphism $\Theta\colon \tE'\rightarrow \tE$ 
which fits into a commutative diagram 
\begin{equation}\label{rft2d}
\xymatrix{
\tE'\ar[d]^{\tau}\ar[rd]^{\sigma'} \ar@/_1pc/[dd]_{\Theta}&\\
\tG\ar[d]^-(0.5){\upgamma}\ar[r]^-(0.4){\pi}\ar@{}|\Box[rd]&X'_\et\ar[d]^{g}\\
\tE\ar[r]^{\sigma}&X_\et}
\end{equation}
{\em We prove that the lower square is Cartesian.} 

We prove first a base change theorem relatively to this square for torsion abelian sheaves of $X'_\et$,
inspired by a base change theorem for oriented products due to Gabber (\cite{ag1} 6.5.5). 
It reduces to the proper base change theorem for the étale topos. Then, we prove the following result which plays a crucial role in the proofs of both \ref{fpscpdi4} and \ref{fpscpdi7}: 

\begin{teo}[\cite{ag1} 6.5.31]\label{rft3} 
Let $g\colon X'\rightarrow X$ be a {\em proper} and smooth morphism of smooth $S$-schemes. 
Then, there exists an integer $N\geq 0$ such that for any integers $n\geq 1$ and $q\geq 0$, 
and any quasi-coherent $\co_{X'_n}$-module $\cF$, the kernel and cokernel of the base change morphism 
\begin{equation}\label{rft3a} 
\sigma^*(\rR^q g_*(\cF))\rightarrow \rR^q\upgamma_*(\pi^*(\cF)),
\end{equation}
are annihilated by $p^N$.
\end{teo}
In this statement, $\pi^*$ and $\sigma^*$ denote the pull-backs for the morphisms of ringed topos 
\begin{eqnarray}
\pi\colon (\tG,\tau_*(\ocB'))&\rightarrow &(X'_\et,\co_{X'}),\label{rft3b}\\
\sigma\colon (\tE,\ocB)&\rightarrow &(X_\et,\co_{X}).\label{rft3c} 
\end{eqnarray}


\begin{thebibliography}{99} 

\bibitem{ag1} {\sc A. Abbes, M. Gros},  Les suites spectrales de Hodge-Tate, preprint (2000), 
\href{https://arxiv.org/abs/2003.04714}{arXiv:2003.04714}.

\bibitem{ag2} {\sc A. Abbes, M. Gros},  Correspondance de Simpson $p$-adique II, Fonctorialité par image directe propre et Systèmes locaux de Hodge-Tate, in preparation. 

\bibitem{agt} {\sc A. Abbes, M. Gros, T. Tsuji},  The $p$-adic Simpson correspondence, 
\href{http://press.princeton.edu/titles/10779.html}{Ann. of Math. Stud., {\bf 193}}, Princeton Univ. Press (2016).  

\bibitem{achinger} {\sc P. Achinger},  {\em $K(\pi,1)$-neighborhoods and comparison theorems}, Compositio Math. {\bf 151} (2015), 1945-1964. 

\bibitem{sga4} {\sc M. Artin, A. Grothendieck, J. L. Verdier}, Théorie des topos et cohomologie étale des schémas, SGA 4, 
Lecture Notes in Math. Tome 1, {\bf 269} (1972); Tome 2, {\bf 270} (1972); Tome 3, {\bf 305} (1973), Springer-Verlag. 

\bibitem{deligne1} {\sc P. Deligne}, {\em Théorème de Lefschetz et critères de dégénérescence de suites spectrales}, 
Pub. Math. IH\'ES {\bf 35} (1968), 107-126.

\bibitem{faltings1} {\sc G. Faltings}, {\em $p$-adic Hodge theory}, J. Amer. Math. Soc. {\bf 1} (1988),  255-299. 

\bibitem{faltings2} {\sc G. Faltings}, {\em Almost étale extensions}, in 
Cohomologies $p$-adiques et applications arithmétiques. II, Astérisque {\bf 279} (2002), 185-270. 

\bibitem{faltings3} {\sc G. Faltings}, {\em A $p$-adic Simpson correspondence},
Adv. Math. {\bf 198} (2005),  847-862. 

\bibitem{ega4} {\sc A. Grothendieck, J.A. Dieudonné}, \'Eléments de Géométrie Algébrique, IV \'Etude locale des schémas et 
des morphismes de schémas, Pub. Math. IH\'ES {\bf 20} (1964), {\bf 24} (1965), {\bf 28} (1966), {\bf 32} (1967).

\bibitem{he1} {\sc T. He}, {\em Cohomological descent for Faltings' $p$-adic Hodge theory and applications}, 
preprint (2021), \href{https://arxiv.org/abs/2104.12645}{arXiv:2104.12645}.

\bibitem{he2} {\sc T. He}, {\em Sen operators and Lie algebras arising from Galois representations over $p$-adic varieties}, 
preprint (2022), \href{https://arxiv.org/abs/2208.07519}{arXiv:2208.07519}.

\bibitem{jannsen} {\sc U. Jannsen}, {\em Continuous étale Cohomology}, Math. Ann. {\bf 280} (1988), 207-245. 

\bibitem{ks2} {\sc M. Kashiwara, P. Schapira}, Categories and sheaves, Grundlehren der mathematischen Wissenschaften {\bf 332}, 
Springer-Verlag (2006). 

\bibitem{katz2} {\sc N. Katz}, {\em  Algebraic solutions of differential equations (p-curvature and the Hodge filtration)},
Invent. math. {\bf 18} (1972), 1-118.

\bibitem{lz}  {\sc R. Liu, X. Zhu}, {\em Rigidity and a Riemann-Hilbert correspondence for $p$-adic local systems}, 
Invent. math. {\bf 207} (2017), 291-343. 

\bibitem{ov} {\sc A. Ogus, V. Vologodsky}, {\em Non abelian Hodge theory in characteristic $p$}, Pub. Math. IHES {\bf 106} (2007), 1-138.

\bibitem{tsuji1} {\sc T. Tsuji}, {\em $p$-adic étale cohomology and crystalline cohomology
in the semi-stable reduction case}, Invent. math. {\bf 137} (1999), 233-411. 

\bibitem{tsuji5} {\sc T. Tsuji}, {\em  \href{https://link.springer.com/article/10.1007/s00208-018-1655-2}{Notes on the local $p$-adic Simpson correspondence}}, 
Math. Annalen (2018) {\bf 371} 795-881. 


\end{thebibliography}
\end{document}